


%


\magnification\magstep1 \baselineskip = 18pt \def \comp{ \;{}^{
{}_\vert }\!\!\!{\rm C}   } \def\N#1{\left\Vert#1\right\Vert}  
\def\tore{{\bf T}} \def\ie{{\it i.e.\/}\ } \def\cf{{\it cf.\/}\ }
\def\n{\noindent} \def\qed{{\hfill{\vrule height7pt width7pt
depth0pt}\par\bigskip}} \def\pf{ \medskip \n {\bf Proof.~~}}
\def\lgar{\longrightarrow} \def\T{{\bf T}} \def\ph{\varphi}

\def\vp{\varepsilon} \def\eps{\vp}

\def\n{\noindent} \def\qed{{\hfill{\vrule height7pt width7pt
depth0pt}\par\bigskip}} \def\pf{ \medskip \n {\bf Proof.~~}}
\def\lgar{\longrightarrow} \def\vp{\varepsilon} \def\eps{\vp}
\def\T{{\tore}} \def\ph{\varphi}

\centerline{\bf A polynomially bounded operator on Hilbert space which
is} \centerline{\bf not similar to a contraction}\bigskip
\centerline{\bf by Gilles Pisier}\bigskip\bigskip

\centerline{ (latest revision 27th of February  1996)}\bigskip
\bigskip\bigskip {\bf Abstract} : Let $\eps >0$. We prove that there
exists an operator $T_\eps:\ell_2\to\ell_2$, such that for any
polynomial $P$ we have $\|{P(T)}\| \leq(1+\eps)\|{P}\|_\infty$, but
which is not similar to a contraction, {\it i.e.} there does not exist
an invertible operator $S:\ \ell_2\to\ell_2$ such that $\|{S^{-1}T_\eps
S}\|\leq 1$. This answers negatively  a question attributed to Halmos
after his well known 1970 paper (``Ten problems in Hilbert space").

\vskip24pt \vskip24pt

{\bf Contents: }

\n {  \S 0. Introduction}

\n {  \S 1. Main results}

\n {  \S 2. Martingales}

\n {  \S 3. Proof of Theorem 1.1}

\n {  \S 4. Operator space interpretations}

\vfill\eject

\n {\bf \S 0. Introduction}

\def\Pi{\tore}

\n  Let $D=\{z\in \comp\ |\ |z|<1\}$ and $\Pi=\partial
D=\{z\ |\ |z|=1\}$.

\n Let $H$
  be a  Hilbert space.

\n Any operator $T:H\to H$ with $\|T\| \le 1$ is called a contraction.
By a celebrated inequality due to von Neumann [vN], we have then, for
any polynomial $P$ $$\|P(T)\|\le \sup_{z\in D} |P(z)|.\leqno(0.1)$$ We
say that $T$ is similar to a contraction if there is an invertible
operator $S:\ H\to H$ such that $S^{-1}TS$ is a contraction.

\n An operator $T:H\to H$ is called power bounded if $\sup_{n\geq
1}\N{T^n}<\infty$.

\n It is called polynomially bounded if there is a constant $C$ such
that, for any polynomial $P$ $$\N{P(T)}\leq C
\sup\{|P(z)|\ |\ z\in\comp,\  |z|=1\}.\leqno(0.1)'$$ Clearly, if $T$ is
similar to a contraction, then it is  power bounded,
 and actually, by von Neumann's  inequality, it is polynomially
bounded.

 We denote $$\N{P}_\infty=\sup_{z\in\partial D}|P=AC(z)|.$$ Let $A$ be
the disc algebra, \ie the closure of the space of all
 (analytic) polynomials in  the space  $C(\Pi)$ of all continuous
 functions on $\Pi$, and let $L^p(\tore)$  be the $L^p$-space relative
to the normalized Lebesgue measure on $\tore$.

\n Let $H^p=\{f\in L^p(\Pi)|\ \hat f(n)=0~~\forall~n<0\}$. If $X$ 
is a
Banach space, we denote, for $p<\infty$,  by $H^p(X)$ the analogous
space of $X$-valued functions.  In the particular case $X=B(H)$ with
$H=\ell_2(I)$ ($I$ being an arbitrary set), we denote by
$L^\infty(B(H))$ the space of all (classes of)  bounded  $B(H)$-valued
functions of which all matrix coefficients are measurable.  This space
can be identified isometrically with the dual of the space
$L^1(B(H)_*)$ of all Bochner integrable functions
 with values in the predual of $B(H)$, \ie the space of all trace class
operators.

If $T$ is polynomially bounded, the map $T\to P(T)$ extends to a unital
homomorphism $u_T:A\to B(H)$ and (0.1)' implies $\N{u_T}\leq C$.

In 1959, B\'ela Sz.-Nagy [SN]asked whether every power bounded
operator $T$ is similar to a  contraction. He proved that the answer is
positive if $T$ is invertible and if both $T$ and its inverse are power
bounded (then $T$ is actually similar to a unitary operator, see
[SN]).  In 1964, S. Foguel [Fo] (see also [Ha2]) gave a  counterexample
to Nagy's question using some properties of Hadamard-lacunary Fourier
series.
    Foguel's example is not  polynomially bounded (see [Le]), whence
    the next question:

 \n Is every  polynomially bounded operator similar to a  contraction
?

\n This revised version of  Sz.-Nagy's original question was
popularized by P. Halmos in [Ha1], and since then, many authors refer
to it as``the Halmos problem".

In [Pa2], Paulsen gave a useful criterion for an operator $T$ to be
similar to a contraction. He proved that this holds iff the
homomorphism $u_T$ is ``completely bounded" (see below for more
background on this notion).  In that case, the operator $T$ is called
``completely polynomially bounded".  In these terms, the above problem
becomes: is every polynomially bounded operator completely polynomially
bounded ?  Or equivalently, is every bounded unital homorphism $\pi:
A\to B(H)$ automatically completely bounded ?

In [Pe1, Pe2] (see also [FW]), V. Peller proposes a candidate for a
counterexample: let
 $\Gamma:H^2\to(H^2)^*$ be a  Hankel operator , \ie such that the
 associated bilinear map (denoted again by $\Gamma$) on $H^2\times H^2$
satisfies $\forall~f\in A,\ \forall g,h\in
H^2~~~~\Gamma(fg,h)=\Gamma(g,fh).$ In other words, if we denote for
$f\in A$, by $M_f:H^2\to H^2$ the operator of multiplication by $f$ and
by $^tM_f:(H^2)^{*}\to (H^2)^{*}$ its  adjoint, we have $$\Gamma
M_f={}^t M_f\Gamma~~~~~~~~~~~~~\forall~f\in A.\leqno(0.2)$$ Let
$H=(H^2)^*\oplus H^2$. For any polynomial $f$   consider the operator
$R(f): H\to H$ defined by the following block matrix:
$$R(f)=\left(\matrix{^tM_f&\Gamma M_{f'}\cr 0& M_f\cr
}\right)\leqno(0.3)$$ Since the coefficient   $f\to\Gamma M_{f'}$
behaves like a  derivation, we easily verify that $f\to R(f)$ is a
unital homomorphism on polynomials,  \ie we have $R(1)=1$, and,   by
(0.2), for all polynomials $f,g$, we have
$$R(fg)=R(f)R(g).\leqno(0.4)$$ Hence, under the condition
$$\exists~C~~~\forall~f\ {\rm polynomial}~~~~~~\N{\Gamma M_{f'}}\leq
C\N{f}_\infty,\leqno(0.5)$$ the mapping $f\to R(f)$ defines a bounded
unital homomorphism on $A$, or equivalently the operator
$T_\Gamma=R(z)$ (\ie $R(\ph_0)$ for $\ph_0$ defined by $\ph_0(z)=z$ for
all $z$) is polynomially bounded.

By a  well known theorem of Nehari (see e.g. [Ni]),  each  Hankel
operator $\Gamma$ is associated to a symbol $\ph\in L^\infty$   so that
$\Gamma=\Gamma_\ph$. (This corresponds to the case $\dim(H)=1$ in (0.6)
below.)  We denote simply by $T_\ph$  the  operator $T_\Gamma$ with
$\Gamma=\Gamma_\ph$. In [Pe2], Peller shows that if
 $\ph'\in BMO$ then we have  (0.5) and consequently  $T_\ph$ is
polynomially bounded. The question whether this implies ``similar to a
contraction" was then posed by Peller [Pe2] , but Bourgain [Bo2] showed
that $\ph'\in BMO$ implies that $T_\ph$ is similar to a  contraction,
and very recently (summer 95), Aleksandrov and Peller [AP] showed that
actually $T_\ph$ is polynomially bounded only if $\ph'\in BMO$. In
conclusion, there is no counterexample in this class.

However, as is well  known (cf. e.g. [Ni]) there is a {\it vectorial}
version  of Hankel operators: given a function  $\ph\in
L^\infty(B(H))$, we can associate to it an operator
 (usually called a vectorial Hankel operator) $\Gamma_\ph:H^2(H)\to
 H^2(H)^*$ still satisfying the identity (0.2), with respect to the
multiplication operator $M_f$ considered as acting on $H^2(H)$.  More
precisely, we can define for any $g,h$ in $H^2(H)$
$$\Gamma_\ph(g,h)=\int [\ph(\xi) g(\xi),h(\xi)] dm(\xi),\leqno(0.6)$$
where $[.,.]$ is a bilinear map on $H$ associated to a fixed isometry
between $H$ and $H^*$.  Thus, defining $R(f)$ as in (0.3) above, we
still have (0.4) and the associated  operator
 $T_\ph:H^2(H)^*\oplus H^2(H)\to H^2(H)^*\oplus H^2(H)$ is polynomially
bounded provided $\Gamma=\Gamma_\ph$ satisfies (0.5).

  The main result of this paper is that, when $H$ is infinite
  dimensional, there are counterexamples to the Halmos problem of the
form $T=T_\ph$, with $\ph\in L^\infty(B(H))$.

\n {\bf Remark.}\ There is   some recent  related work (on  operators
of the form $T_\Gamma$) by  Stafney [St] and also by S. Petrovic,  V.
Paulsen and Sarah  Ferguson whose work I have  heard of, through the
seminar talks they gave in Texas at various occasions in 95.

\n {\bf Notation and background.} In general, any unexplained notation
is standard. Let $H,K$ be two Hilbert spaces. We denote by $B(H,K)$
(resp. $B(H)$) the space of all bounded linear operators from $H$ to
$K$ (resp. from $H$ to $H$)  equipped with its usual norm. We denote by
$H^*$ the dual of $H$ which, of course, is a Hilbert space canonically
identifiable with the complex conjugate $\bar H$ of $H$.

\n We denote by $m$ the normalized Lebesgue measure on the unit circle
$\T$.

\n Recall that we denote simply by $A$ the disc algebra.  Note that $A$
is a closed subalgebra of $H^{\infty}$. Often, we implicitly consider
 a function in  $H^\infty$ as extended analytically inside the unit
 disc, in such a way that we recover the original function on the
circle by taking radial (or non-tangential) limits almost everywhere.
When the original function on the circle is actually in $A$, its
analytic extension is continuous on $\bar D$.

\n Consider a function $f$ in (say) $L^p(\T,m)$ ($1\le p\le \infty$).
When we sometimes abusively say that $f$ is ``analytic", what we really
mean is that $f\in H^p$, \ie that $f$ extends analytically inside $D$.
In that case, throughout this paper, the derivative $f'$ of $f$ always
means the derivative of the Taylor series of $f$.

We denote by $\ell_2^n$ the $n$-dimensional Hilbert space, and by
$\ell_2^n(H)$ the (Hilbertian) direct sum of $n$ copies of $H$.
Moreover, we denote $M_n=B(\ell_2^n)$.

Let $H,K$ be Hilbert spaces.  Let $S\subset
 B(H)$ be a subspace.  In the theory of operator algebras, the notion
 of complete boundedness for a linear map $u:S\rightarrow B(K)$ has
been extensively
 studied recently. Its origin lies in the work of Stinespring (1955)
   and Arveson (1969) on completely positive maps (see [Pa1, Pi4] for
more details
 and references). Let us equip $M_n(S)$ and $M_n(B(K))$ (the spaces of
matrices with entries respectively in $S$ and $B(K)$) with the norm
induced respectively by $B(\ell_2^n(H))$ and $B(\ell_2^n(K)).$

\n  A map $u\colon S\to B(K)$ is called completely bounded (in short
c.b.) if there is a constant $C$ such that the maps $I_{M_n}\otimes u$
are uniformly bounded by $C$ i.e. if we have $$\sup_n \| I_{M_n}\otimes
u\|_{M_n(S)\to M_n(B(K))}\le C,$$ and the c.b. norm $\|u\|_{cb}$ is
defined as the smallest constant $C$ for which this holds.

\n When $\|u\|_{cb}\le 1$, we say that $u$ is completely contractive
(or a complete contraction).

\n Let ${\cal A}\subset B(H)$ be a unital subalgebra, and let $\pi:
{\cal A}\to B(H)$ be a unital homomorphism. Paulsen ([Pa2]) proved that
$\pi$ is completely bounded iff there is an invertible operator $S:
H\to H$ such that $S^{-1}\pi(.)S$ is completely contractive.

 We now return to polynomially bounded operators. An operator $T\colon
 \ H\to H$ will be called completely polynomially bounded  if there is
 a constant $C$ such that for all $n$ and all  $n\times n$ matrices
 $(P_{ij})$ with polynomial entries we have
 $$\|(P_{ij}(T))\|_{B(\ell^n_2(H))} \le C\sup_{z\in
 \T}\|(P_{ij}(z))\|_{M_n}\leqno (0.7)$$ where $(P_{ij}(T))$ is
 identified with an operator on $\ell^n_2(H)$ in the natural way.
 Recall that $M_n$ is identified with $B(\ell_2^n)$.
 Note that $T$ is completely polynomially bounded  iff the homomorphism
 $P\to P(T)$ defines a completely bounded homomorphism $u_{T}$ from the
 disc algebra $A$ into
 $B(H)$. Here of course we consider $A$ as a subalgebra of the
 $C^*$-algebra $C({\bf T})$ which itself can be embedded e.g. in
 $B(L_2({\bf T}))$ by identifying a function $f$ in $C({\bf T})$ or
 $L_\infty({\bf T})$ with the operator of multiplication by $f$ on
 $L_2({\bf T})$.

 We can now state Paulsen's criterion :

 \proclaim Theorem 0.1. ([Pa2]) An operator $T$ in $B(H)$ is similar to
a
 contraction iff it is completely polynomially bounded. Moreover $T$ is
 completely polynomially bounded  with constant $C$ (as in (0.7) above)
 iff there is an isomorphism
 $S\colon \ H\to H$ such that $\|S\|\, \|S^{-1}\|\le C$ and
 $\|S^{-1}TS\| \le 1$.

Actually, we only use the easy direction of this criterion, which can
be derived  as follows from Sz.-Nagy's well known unitary dilation
theorem.  Let $T: H\to H$ be a contraction. Sz.-Nagy proved (see [SNF])
that we can find a larger Hilbert space  $\hat H$ containing $H$ as a
subspace, and a unitary  operator $U$ on  $\hat H$, such that $$\forall
n\ge 0\quad T^n=P_H U^n _{|H}.$$ (Here $P_H$ denotes the orthogonal
projection from  $\hat H$ to $H$.) Hence, we have $P(T)=P_H P(U)
_{|H}$, for all polynomials $P$. In particular, this implies von
Neumann's inequality (0.1).  More generally,  for any matrix $[P_{ij}]$
with polynomial entries we have $[P_{ij}(T)]= j^*[P_{ij}(U)]j$ where
$j^*$ is the diagonal matrix with diagonal entries all equal to $P_H$.
Therefore, $$\|[P_{ij}(T)]\|\le \|[P_{ij}(U)]\| \le \sup_{z\in
\T}\|(P_{ij}(z))\|_{M_n}.$$ The last inequality is a direct consequence
of the fact that $U$ generates  a commutative $C^*$-subalgebra of
$B(\hat H)$, and its spectrum lies in $\T$.  Clearly, this implies
$$\|[P_{ij}(S^{-1}TS)]\| \le \|S\|\|S^{-1}\|\sup_{z\in
\T}\|(P_{ij}(z))\|_{M_n},$$ which proves that similarity to a
contraction implies complete polynomial boundedness.

For any linear map $u: A\to B(H)$ it is easy to check that, for any $n$
and
 for any finitely supported sequence $(a_k)$  in $M_n$, we can write:
$$\left\|\sum a_k\otimes u(z^k)\right\|_{B(\ell_2^n(H))}\le \|u\|_{cb}
\sup_{z\in\T} \left\|\sum a_k\ z^k\right\|_{M_n}.\leqno(0.8)$$ This
formulation   will be used below.

We refer the reader to the books [SNF, Pa1] for more background on
dilation Theory, and to  the  lecture notes volume [Pi4] which
describes the ``state of the art" on similarity problems until 95.

\vskip24pt \n {\bf \S 1. Main results}

\n Let $H$ be a Hilbert space. Consider an  operator $$\Gamma\colon
\ H^2(H)\longrightarrow H^2(H)^*,$$ or equivalently  a bounded bilinear
form $\Gamma\colon \ H^2(H) \times H^2(H) \longrightarrow {\bf C}.$
Then $\Gamma$ is called Hankelian (or a  Hankel operator) if for any
multiplication operator $M_\varphi\colon \ H^2(H)\lgar H^2(H)$ by a
polynomial or a function $\varphi$ in $A$ (or in $H^\infty$), we have
$$\forall g,h\in H^2(H)\qquad\Gamma(g\varphi ,h) = \Gamma(g,\varphi h).
\leqno (1.1)$$ \n Equivalently we have $\Gamma M_\ph={^t M_\ph\Gamma}.$
Let $F\to{\cal D}(F)$ be a linear mapping from $A$ into
$B(H^2(H),H^2(H)^*)$ which is a ``derivation" in the following sense :
$$\forall F,G \in A\qquad {\cal D}(FG)={^tM_F}{\cal D}(G)+{\cal D}(F)
M_G.\leqno(1.2)$$ Let ${\cal H}=H^2(H)^*\oplus H^2(H)$. Then the
mapping $R:A\to B({\cal H})$ defined by
$$R(F)=\left(\matrix{{^tM_F}&{\cal D}(F)\cr 0& M_F\cr}\right)$$ clearly
is  a unital homomorphism.

Now assume in addition that ${\cal D}(F)$ is Hankelian for any $F$ in
$A$, \ie that ${\cal D}(F)M_\ph={^tM_\ph}{\cal D}(F)$ for all $\ph$ in
$A$. Then a simple computation shows by induction that $${\cal
D}(F^2)=2{\cal D}(F)M_F~{\rm and}~{\cal D}(F^n)=n{\cal D}(F)
M_{F^{n-1}}$$ for all $n\geq 1$. Therefore, for any polynomial $P$ we
must have $${\cal D}(P(F))={\cal D}(F)M_{P'(F)},$$ where $P'$ is the
derived polynomial.  Applying this with the function $F=\ph_0$ defined
by $\ph_0(z)=z$ for all $z$, we find $${\cal D}(P)={\cal
D}(\ph_0)M_{P'}.$$ So if we let $\Gamma={\cal D}(\ph_0)$ we have for
all polynomials $P$ $${\cal D}(P)=\Gamma M_{P'}.\leqno(1.3)$$
Conversely, given any Hankel operator $\Gamma$, if we define ${\cal
D}_{\Gamma}$ by setting ${\cal D}_{\Gamma}(P)=\Gamma M_{P'}$, then
${\cal D}_{\Gamma}(P)$ is Hankelian, satisfies the derivation identity
(1.2) for all polynomials $F,G$, and ${\cal D}_{\Gamma}(\ph_0)=\Gamma$.
(Thus (1.3) defines a one to one correspondence between $\Gamma$ and
${\cal D}_{\Gamma}$.)

\n Given    $\Gamma$   (and the associated ${\cal D}_{\Gamma}$), let
$T_\Gamma\colon \ H^2(H)^* \oplus H^2(H) \lgar H^2(H)^* \oplus H^2(H)$
be the operator defined by the operator matrix $$T_\Gamma =
\left(\matrix{{}^t S&\Gamma\cr 0&S\cr}\right)$$ where $S\colon
\ H^2(H)\lgar H^2(H)$ is the shift operator i.e.\ $S = M_{\varphi_0}$,
and ${}^t S$ denotes its adjoint on the dual space $H^2(H)^*$.  We will
show in this paper that, in the {\it vectorial\/} case (with
$\dim(H)=\infty$), there are p.b.\ operators of this type which are not
similar to a contraction. The preceding remarks show that, for any
polynomial $P$ $$P(T_\Gamma) = \left(\matrix{{}^tP(S)&\Gamma M_{P'}\cr
0 &P(S)\cr }\right).$$ Therefore if there is a constant $C$ such that
$$\|\Gamma M_{P'}\| \le C\|P\|_\infty\leqno (1.4)\qquad \forall P$$
then $T_\Gamma$ is polynomially bounded. Our main result is the
following:

\proclaim Theorem 1.1. Let $(C_n)$ be any sequence in $B(H)$ such that
$$\left\|\sum \alpha_nC_n\right\|_{B(H)} \le \left(\sum
|\alpha_n|^2\right)^{1/2},\leqno (1.5)\qquad \forall (\alpha_n)
\in\ell_2$$ and let $(K_n)$ be an increasing sequence of positive
integers such that for all $n> 1$ $$2^{n-1}<K_n \le 2^n.\leqno(1.6)$$
Consider the operator $u\colon \ A\longrightarrow B(H)$ defined by
$u(F) = \sum_{n\ge 1} \hat F(K_n)C_n.\qquad$ \break Then there is a
function $\ph\in L^\infty(B(H))$ such that the associated vectorial
Hankel operator $\Gamma=\Gamma_\ph \colon \ H^2(H) \lgar H^2(H)^*$
 satisfies (0.5) for some constant $C$ and moreover is such that for
any polynomial $P$ (with derived polynomial denoted by  $P'$):  $$
\forall x,y\in H \quad \langle u(P) x,\overline y\rangle = \Gamma
[{P'}(1\otimes x)] (1\otimes y) \leqno(1.7) $$
 where we denote by $1\otimes x$ the element of $H^2(H)$ corresponding
to the function taking constantly the value $x$, and where $y\to
\overline y$ is an antilinear isometry on $H$.

\proclaim Corollary 1.2. There is a polynomially bounded operator (of
the form $T_\Gamma$) which is not similar to a contraction.

\pf Let $(C_n)$ be a sequence of operators satisfying the CAR
(``canonical anticommutation relations''), i.e.\ such that $$C_iC_j +
C_jC_i=0\quad C^{*}_iC_j + C_jC^{*}_i = \delta_{ij}I.  \leqno \forall
i,j$$ It is well known that this implies (1.5) (with equality even),
see e.g. [BR, p.  11]. Moreover, if $(C'_1,...,C'_n)$ is another
$n$-tuple satisying the CAR, there is a $C^*$-representation taking
$(C_1,...,C_n)$ to $(C'_1,...,C'_n)$ (see [BR, Th. 5.2.5]), so that for
any $a_1,...,a_n$ in $M_N$ (with $N\ge 1$) the value of $\|\sum_1^n
a_i\otimes C_i\|_{B(\ell^N_2(H)}$ does not depend on the particular
realization of $(C_1,...,C_n)$. We denote by $\bar x$ the complex
conjugate of a matrix $x$. Equivalently, if $x\in B(H)$, we can
identify $\bar x$ with ${}^t x^*: H^*\to H^*$.

\n  We claim that, with this choice of $(C_n)$,  the mapping $u$ is
{\it not\/} completely bounded.  To see this last fact (also well
known, see e.g. [H, example 3.3]), we first fix $n$ and recall that
there is (using the Pauli or Clifford  matrices, see also [BR, p. 15])
a  realization of $(C_1,\ldots, C_n)$ satisfying the CAR in the space
of (say) $2^n\times 2^n$ matrices.  Then observe that we have by (0.8)
$${\left\|\sum^n_1 C_k \otimes \overline C_k\right\| = \left\|\sum
u(z^{K_k}) \otimes \overline C_k\right\| \le \|u\|_{cb} \sup_{|z|=1}
\left\|\sum^n_1 z^{K_k} \otimes \overline C_k \right\| = \|u\|_{cb}
\sqrt n.}$$ However, since we can assume $(C_1,\ldots, C_n)$ realized
in the space of  $2^n\times 2^n$ matrices, of which the unit matrix is
denoted by $I$, we have $$\left\|\sum\nolimits^n_1 C_k\otimes \overline
C_k\right\| \ge {\hbox{tr} \left(\sum\nolimits^n_1
C_kC^{*}_k\right)\over \hbox{tr } I} = {n\over 2}.$$ (To check this,
note  that $$\left\|\sum\nolimits^n_1 C_k\otimes \overline C_k\right\|
= \sup \left|\sum\nolimits^n_1 \hbox{ tr } (C_k X C^*_k Y)\right|$$
where the supremum runs over all $2^n\times 2^n$ matrices $X,Y$ with
Hilbert-Schmidt norm 1 and take $X=Y=I(\hbox{tr } I)^{-1/2}$ .) Thus 
we
conclude that $\|u\|_{cb} \ge \sqrt n/2$ for all $n$, hence $u$ is {\it
not\/} completely bounded.  Therefore, if $\Gamma$ satisfies (1.7), the
mapping $F\to \Gamma M_{F'}$ cannot be c.b., hence a fortiori $F\to
F(T_\Gamma)$ is not c.b. on $A$, which ensures by Paulsen's criterion
(easy direction) [Pa2] that $T_\Gamma$ is not similar to a
contraction.\qed

\n {\bf Remark 1.3.} Let $T_\Gamma = \left(\matrix{{}^tS&\Gamma\cr
0&S\cr}\right)$ be an operator as in the last corollary. Fix $\vp>0$,
and consider $$T_{\vp\Gamma}  = \left(\matrix{{}^t S&\vp\Gamma\cr
0&S\cr}\right).$$ Then, by (1.4), we have for all polynomials $P$
$$\|P(T_{\vp\Gamma})\| \le (1+\vp C) \|P\|_\infty,$$ hence we can get
an operator with polynomially bounded constant arbitrarily close to 1,
but still not similar to a contraction.

\n {\bf Remark 1.4.} In the preceding argument for Corollary 1.2, we
are implicitly using the ideas behind a joint result of V.~Paulsen and
the author (included in [Pa3, Th. 4.1]). The latter result shows in
particular the following:\ let $m = (m(n))_{n\ge 0}$ be a scalar
sequence and let $u_m\colon \ A \lgar B(H)$ be the operator defined by
$u_m(F) = \sum\limits_{n\ge 0} \hat F(n) C_n m(n)$, with $(C_n)$
satisfying the CAR. Then $u_m$ is c.b.\ iff $\sum |m(n)|^2<\infty$.

Let $(m(n))_{n\geq 0}$ be a scalar sequence such that $$\sup_{n\geq
0}\sum_{2^{n-1}<k\leq 2^n}|m(k)|^2<\infty.\leqno(1.8)$$ It is well
known that this condition characterizes the Fourier multipliers from
$H^1$ to $H^2$ on the circle (\cf e.g. [D, p. 103]).

\n With essentially the same arguments as for Theorem 1.1, we can prove
more generally

\proclaim Theorem 1.5. Let $(C_n)$ be as in Theorem 1.1. Assume that
$(m(n))_{n\geq 0}$ satisfies (1.8) as above. Let $u: A\to B(H)$ be the
mapping defined by $$u(F)=\sum_{n\geq 0}\hat F(n)m(n) C_n.\leqno(1.9)$$
Then the conclusion of Theorem 1.1 still holds.

\n See Remark 3.3 below for an indication of proof.

\n {\bf Remark 1.6.} Let $\ph\in L^\infty(B(H))$ be the
 function associated to (1.9) by Theorem 1.5.  By (0.6) and (1.7) the
Fourier transform of $\ph$ restricted to non positive integers is
determined and we have for all $n\ge 1$ $$C_n m(n)= \int \ph(e^{it}) n
e^{i(n-1)t} dm (e^{it}).$$ Hence as a ``symbol" in the sense of
Nehari's theorem as in (0.6), we can take just as well
$$\ph(e^{it})=\sum_{n\ge 1} n^{-1} C_{n} m(n) e^{-i(n-1)t}.$$ Let
$f=\sum_{k\ge 0} \hat f(k)z^k \in A$. Let $P_n(f) =\sum_{0\le k\le n}
\hat f(k)z^k.$ Let us denote by $(e_n)$ the canonical basis in
$\ell_2$.  Then,  assuming (1.8)
 the key inequality (0.5) is equivalent to the following one, perhaps
of independent interest: there is a constant $C$ such that for any
analytic  function $f$ in $BMO$ $$\left\|\sum _{n\ge 1} e_{n} n^{-1}
m(n) e^{-i(n-1)t}P_{n-1}(f')\right\|_{BMO(\ell_2)} \le
C\|f\|_{BMO},\leqno(1.10)$$ where we have denoted by $BMO$ (resp.
$BMO(\ell_2)$) the classical space of $BMO$ functions on the circle
with values in $\comp$   (resp. $\ell_2$).

\n Note that (1.10) is actually proved below using martingale versions
of $H^1$ and $BMO$, but the end result can be formulated without
reference to martingales as we just did in (1.10).

\n{\bf Remark 1.7}. Let us examine the converse to Theorem 1.5. Fix a
sequence $(C_n)$ for which there is a constant $\delta>0$ such that
$$\forall~\alpha=(\alpha_n)\in\ell_2~~~~~~\delta(\sum|\alpha_n|^2)^{
1/2}\leq\N{\sum\alpha_n
C_n}.\leqno(1.11)$$ Let $(m(n))_{n\geq 0}$ be any scalar sequence.
Then if the conclusion of Theorem 1.1 holds for the mapping $u$ defined
in (1.9), we have necessarily (1.8).
 This follows from [AP].  Indeed, let $\ph\in L^\infty(B(H))$ be a
function associated to such a $u$ as in Theorem 1.1, with associated
Hankel operator satisfying (0.5).  For any sequence $(\beta_n)$ in the
unit ball of $\ell_2$, we have a linear form $\xi$ in $B(H)^*$ with
$\|\xi\|_{B(H)^*}\le 1/\delta$ such that $\xi(C_n)=\beta_n$, for all
$n$.  Clearly, the function $\xi(\ph)\in L^\infty$ defines a Hankel
operator $\Gamma_{\xi(\ph)}$ still satisfying  (0.5).  By [AP], this
implies that $\xi(\ph')$ is in BMO, or equivalently that $\sum_{n\ge 0}
\beta_n m(n) z^n$ is in BMO. Since this holds for all $(\beta_n)$ in
$\ell_2$, this means, after transposition, that $(m(n))$ is a (bounded)
multiplier from $H^1$ to $\ell_2\approx H^2$. As already mentioned,
this is equivalent to (1.8) (\cf e.g. [D, p. 103]).  In particular, a
sequence $(K_n)$  satisfies the conclusion of Theorem 1.1 (for any
$(C_n)$ satisfying (1.5) and (1.11)) if and only if it is the union of
finitely many sequences which are subsequences of a sequence verifying
(1.6). Equivalently, iff it is a finite union of Hadamard lacunary
sequences.  \qed

Of course, our results can also be described in terms of Hankel and
Toeplitz matrices.

\n Let $G$ be the Hankel matrix (with operator entries) defined by
setting for all $i,j\geq 0$
$$G_{ij}=m(i+j+1)(i+j+1)^{-1}C_{i+j+1}.\leqno(1.12)   $$ For $f\in A$,
let us denote by $T(f)$ the Toeplitz matrix (with scalar entries)
defined by setting for all $i,j\geq 0$ $$T(f)_{ij}=\hat f(i-j).$$ Then
Theorem 1.5 can be reformulated as follows : assuming (1.8), there is a
constant $C$ such that for all $f$ in $A$ we have $$\N{GT(f')}\leq
C\N{f}_A.\leqno(1.13)$$ Note that, actually, the BMO norm of $f$ can be
substituted to $\N{f}_A$ in (1.13), see Remark 3.2 below.

\n To verify (1.12), we simply apply (1.7) to $F(z)=z^n$, to show that
the matrix coefficients of $\Gamma$, with respect to the decomposition
$H^2(H)=\oplus_{i\geq 0}z^i H$ are given by (1.12). Then (1.13) reduces
to (1.4).

\vskip24pt \n {\bf \S 2. Martingales}

We will first prove an analogue of Theorem~1.1 with $A$ replaced by its
martingale version denoted by $A_{\rm m}$. Then, in the next section we
will deduce Theorem 1.1 from the martingale case.

Consider   $\Omega = {\bf T}^I$ with $I = \{1,2,3,\ldots\}$. Let
$(Z_i)_{i\ge 1}$ denote the coordinates on $\Omega$ and let ${\cal
A}_n$ be the $\sigma$-algebra generated by $(Z_1,\ldots, Z_n)$ with
${\cal A}_0$ the trivial $\sigma$-algebra. Every continuous function
$f\colon \ \Omega\to {\bf C}$ defines a martingale $(f_n)_n$ by setting
$f_n = {\bf E}(f\mid {\cal A}_n)$.

A martingale $(f_n)_n$, relative to the filtration $({\cal A}_n)$, is
called ``Hardy'' if for each $n\ge 1$ the function $f_n$ depends
analytically on $Z_n$ (but arbitrarily on $Z_1,\ldots, Z_{n-1}$).

Let $1\le p\le \infty$. We denote by $A_{\rm m}$ (resp. $H^{p}_{\rm
m}$) the subspace of $C(\Omega)$ (resp. $L^{p}(\Omega,P)$) formed by
all $F$ which generate a Hardy martingale. As is well known $A_{\rm
m}$  (resp. $H^{\infty}_{\rm m}$) is a uniform algebra in $C(\Omega)$
(resp. $L^{\infty}(\Omega,P)$). Moreover, for all $F,G$ in
$H^{\infty}_{\rm m}$ we have $$(FG)_n = F_nG_n\qquad \forall n\ge
0.\leqno (2.1)$$ In Harmonic Analysis terms, the space $A_{\rm m}$
(resp. $H^{p}_{\rm m}$) is indeed the version of the disc algebra
(resp. $H^{p}$) associated to the ordered group ${\bf Z}^{(I)}$ (formed
of all the finitely supported families $n=(n_i)_{i\in I}$ with $n_i\in
{\bf Z}$), ordered by the {\it lexicographic\/} order, \ie the order
defined by setting $n'<n''$ iff  the last differing coordinate
(=``letter" with reversed alphabetical order) satisfies $n'_i<n''_i$.
 As explained e.g.\ in [Ru, Chapter 8] or in [HL], this group has a
``{\it linear\/}'' behaviour and the associated $H^p$ spaces on it
behave like the classical (unidimensional) ones.     Let $P$ be the
usual probability measure on ${\bf T}^I$ (= normalized Haar measure).
For any Banach space $X$, we will denote by $H^p_{\rm m}(X)$ $(1\le
p\le \infty)$ the usual $H^p$-space of $X$-valued functions on the
ordered group $\Omega$. Equivalently, for $p<\infty$ $H^p_{\rm m}(X)$
is the closure of $A_{\rm m}\otimes X$ in the space $L_p(\Omega,P;X)$
in Bochner's sense.

Our main goal in this section is the following:  \proclaim Theorem 2.1.
There is a bounded unital homomorphism $$\pi\colon \ A_{\rm m}\lgar
B(H)$$ which is not completely bounded.

For each $F$ in $H^{\infty}_{\rm m}$ we denote by $(F_n)$ the
associated Hardy martingale and we set $dF_n = F_n-F_{n-1}$ $\forall
n\ge 1$, $dF_0 = F_0$.

Let $\eta=(\eta_n)$ be an adapted sequence (\ie $\eta_n$ is ${\cal
A}_n$-measurable for all $n\ge 0$) of bounded random  variables such
that there is a constant $C'$ for which:  $$\forall n\ge 0\quad
\|\eta_n\|_\infty\le C'.\leqno(2.2)$$
 Consider $F\in H^{\infty}_{\rm m}$. We introduce a mapping $${\cal
D}^{\eta}(F)\colon \ H^2_{\rm m}(H)\lgar H^2_{\rm m}(H)^*$$ which,
viewed as a bilinear map on $H^2_{\rm m}(H) \times H^2_{\rm m}(H)$, can
be written as   $\forall g,h\in H^{\infty}_{\rm m}(H)$ $${\cal
D}^{\eta}(F)(g,h)  =  {\bf E} \left( \sum_{n\ge 1}\bar Z_n\,
\eta_{n-1}\, dF_n [C_ng,h] \right)\leqno (2.3)$$ where, for convenience
of notation, we use a {\it bilinear} pairing on $H\times H$ $(x,y)\to
[x,y]$ on the right side of this definition, so that $Z \to [C_ng(Z),
h(Z)]$ is in $H^1$ when $g,h$ are in $H^2_{\rm m}(H)$. We assume that
the associated map $J\colon \ H\to H^*$ is isometric.

\n Concerning the issue of the convergence of the series in (2.3), note
that \break $\|F\|^2_{2}= \sum_{n\ge0}{\bf E}|dF_n|^2$, and we take the
precaution to assume that  $g,h$ belong to $H^{\infty}_{\rm m}(H)$
(which is dense in $H^{2}_{\rm m}(H)$).  Therefore,
 by (1.5) the series $$\sum_{n\ge 1}\bar Z_n\, \eta_{n-1}\, dF_n
[C_ng,h]$$ is convergent in $L^2(\Omega,P)$, so that (2.3) is well
defined.

\n But actually, we will show (see below) that (2.3) defines a bounded
operator from $H^2_{\rm m}(H)$ to its dual, so that (2.3) eventually
will make sense for all $g,h$ in $H^2_{\rm m}(H)$.

We denote by $M_{\rm m}(F)\colon\ H^2_{\rm m}(H)\to H^2_{\rm m}(H)$ the
operator of multiplication by $F$ when $F$ is in $H^{\infty}_{\rm m}$.

The following identity has played a crucial role in guiding us to the
present paper. We refer the reader to [LPP, Pi3] for various results
related to this ``derivation identity''.  $${\cal D}^{\eta}(FG) = {\cal
D}^{\eta}(F) M_{\rm m}(G) + {}^tM_{\rm m}(F){\cal D}^{\eta}(G).\leqno
(2.4)\qquad \forall F,G\in H^{\infty}_{\rm m}$$

Let $u^{\eta}\colon \ H^{\infty}_{\rm m}\lgar B(H)$ be the mapping
defined by $$u^{\eta}(F) = \sum_{n\ge 1} {\bf E}(dF_n\bar Z_n
\eta_{n-1}) C_n$$ where $(C_n)$ is a CAR sequence as before.
 The same argument as in the preceding section shows that, if $\eta$ is
 suitably chosen, for instance if $\eta_n=1$ identically for all $n$,
$u^{\eta}$ is {\it not\/} c.b., even when restricted to $A_{\rm m}$.
Therefore, in this case,
 any homomorphism $\pi$ admitting $u^{\eta}$ as a ``coefficient'' a
fortiori is not c.b.\ either.

Moreover we will see that if $V\colon\ H\lgar H^2_{\rm m}(H)$ is the
mapping defined by $V(x) = 1\otimes x$, then we have (note that
${}^tV\colon \ H^2_{\rm m}(H)^*\to H^*$) $$Ju^{\eta}(F) = {}^tV{\cal
D}^{\eta}(F)V\leqno (2.5)\quad \forall F\in H^{\infty}_{\rm m}$$ where
$J\colon \ H\to H^*$ is the isometric isometry associated to the
bilinear pairing $(x,y)\to [x,y]$ on ${ H}$. From (2.4) and (2.5),
Theorem~2.1 is immediate. Indeed, (2.4) shows that the mapping
$$\pi^{\eta}(F) = \left(\matrix{{}^tM_{\rm m}(F)&{\cal D}^{\eta}(F)\cr
0&M_{\rm m}(F)\cr}\right)$$ is a unital homomorphism. By (2.5), ${\cal
D}^{\eta}$ is not c.b.\ if $\eta$ is suitably chosen,  hence a fortiori
$\pi^{\eta}$ is not c.b.\ and we are done.

\n  The proof of (2.5) is immediate, so it remains, to complete the
proof, to check (2.4) and prove that $F\to {\cal D}^{\eta}(F)$ is
bounded. The key for (2.4) is the observation that for any ${\cal
A}_{n-1}$-measurable function $\varphi_{n-1}$ and for any pair $F,G$
say in $H^{\infty}_{\rm m}$ we have $\forall n\ge 1$ $${\bf E}(\bar Z_n
\, dF_n\, dG_n\, \varphi_{n-1}) = 0.\leqno (2.6)$$ Indeed, integrating
first in $Z_n$ alone we find 0 since $F_n,G_n$ are analytic in $Z_n$,
and $F_{n-1}$ can be obtained from $F_n$ by setting $Z_n=0$ in $F_n$
(Jensen's formula). Now, using (2.6) repeatedly it is easy to verify
that for all $F$ in $H^{\infty}_{\rm m}$ we have $${\cal
D}^{\eta}(F)(g,h) = E\left(\sum_{n\ge 1} \bar Z_n \, \eta_{n-1} \, dF_n
[C_ng_{n-1}, h_{n-1}]\right).\leqno (2.7)$$ Hence, using $$d(FG)_n =
F_nG_n - F_{n-1}G_{n-1} = dF_n\, G_{n-1} + F_{n-1}\, dG_n + dF_n\,
dG_n$$ and using (2.6) and (2.7) again, we find $$\eqalign{{\cal
D}^{\eta}(FG)(g,h) &= E\left(\sum_{n\ge1} \bar Z_n\, \eta_{n-1} \,dF_n
[C_nG_{n-1}g_{n-1}, h_{n-1}]\right)\cr &\quad + E\left(\sum_{n\ge 1}
\bar Z_n\, \eta_{n-1} \,dG_n[C_ng_{n-1}, F_{n-1}h_{n-1}] \right)\cr &=
{\cal D}^{\eta}(F) (Gg,h) + {\cal D}^{\eta}(G)(g,Fh).}$$ This
establishes (2.4). We now turn to the boundedness of $F\to {\cal
D}^{\eta}(F)$ on $H^{\infty}_{\rm m}$.

\n Recall that the predual of $B(H)$ can be naturally identified with
the projective tensor product $H\widehat \otimes H$. Therefore,  if
$g,h$ are in the unit ball of $H^2_{\rm m}(H)$ then $g\otimes h$ can be
viewed as an element of the unit ball of $H^1_{\rm m}(H\widehat\otimes
H)$.  Thus, we can use the easy direction of the ``vectorial Nehari
Theorem"  to argue that $$\|{\cal D}^{\eta}(F)\| \le \left\|\sum \bar
Z_n \, \eta_{n-1} \,dF_n\, C_n\right\|_{H^1_{\rm m}(H\widehat \otimes
H)^*}. \leqno (2.8)$$ Now let $v\colon \ \ell_2\to B(H)$ be the mapping
defined by  $v(e_n) = C_n$ which by (1.5) satisfies $\|v\|\le 1$. Then
(2.8) can be estimated through $\ell_2$:  $$\|{\cal D}^{\eta}(F)\|\le
\left\|\sum \bar Z_n\, \eta_{n-1} \,dF_n \, e_n\right\|_{H^1_{\rm
m}(\ell_2)^*}.\leqno (2.9)$$ On the other hand, let  $H^1_{\rm
M}(\ell_2)$ be the space formed of all $({\cal A}_n)$-adapted
$\ell_2$-valued martingales $f = (f_n)_{n\ge 0}$ such that $f^* =
\sup\|f_n\|$ is integrable.  We equip $H^1_{\rm M}(\ell_2)$ with the
norm $f\to {\bf E}(f^*)$. We claim that the natural  inclusion
$$H^1_{\rm m}(\ell_2)\lgar H^1_{\rm M}(\ell_2)$$ taking $f$ to the
associated martingale $(f_n)$ is bounded. Indeed, this is    entirely
classical: let  $f\in H^1_{\rm m}(\ell_2)$.   Since $(f_n)$ is
``analytic" (= a Hardy martingale), by Jensen's inequality  the
sequence $(\|f_n\|^p)$ is a submartingale for any $p>0$, then choosing
$p=1/2$, the desired result follows from Doob's maximal inequality in
$L_2$ (cf. e.g. [Du, p. 307]). This proves the claim.
 Thus, $f\to {\bf E}(f^*)$ is an equivalent norm on $H^1_{\rm
m}(\ell_2)$. As is well known ([Bu, Du, Ga]), the latter norm is
equivalent to the norm $f \to {\bf E}\left(\sum\limits_{n\ge
0}\|df_n\|^2\right)^{1/2}$.

\n Moreover, the dual of $H^1_{\rm M}(\ell_2)$ can be identified with
the space $\hbox{BMO}_{\rm M}(\ell_2)$ defined as the space of all
martingales $y = (y_n)_{n\ge 0}$ such that for all ${n\ge 1}$
$$\sup_{n\ge 1} \left\|{\bf E}(\|y-y_{n-1}\|^2 \mid {\cal
A}_n)\right\|_\infty <\infty$$ equipped with the ``norm" (up to an
additive constant) $$|||y||| = \left(\sup_{n\ge 1} \|{\bf E}
(\|y-y_{n-1}\|^2\mid {\cal A}_n)\|_\infty\right)^{1/2} .$$ It follows
from these remarks that there are numerical constants $K'$ and $ K$
such that if we denote by $(e_n)$ the canonical basis of $\ell_2$ and
if we set $y = \sum \bar Z_n\, \eta_{n-1} \,dF_n\, e_n$ then we have
$$\leqalignno{\|y\|_{H^1_{\rm m}(\ell_2)^*} &\le K'\|y\|_{H^1_{\rm
M}(\ell_2)^*} &(2.10)\cr &\le K[ |||y|||+({\bf E}\|y\|^2)^{1/2}]}$$ But
by a simple computation   we find (recalling (2.2) and denoting $E_n$
for $E(~~~\mid {\cal A}_n)$) $$\eqalign{{\bf E}_n \|y-y_{n-1}\|^2 &=
{\bf E}_n\left\|\sum_{k\ge n} e_k(\bar Z_k \, \eta_{k-1} \,dF_k -
E_{n-1}(\bar Z_k \, \eta_{k-1} \,dF_k))\right\|^2\cr &\le 2{C'}^2{\bf
E}_n\left(\sum_{k\ge n} |dF_k|^2 + E_{n-1} \sum_{k\ge n}
|dF_k|^2\right)\cr &\le 2{C'}^2\left( E_n|F-F_{n-1}|^2
+E_{n-1}|F-F_{n-1}|^2\right)\cr &\le 4{C'}^2|||F|||^2\cr &\le
16{C'}^2\|F\|^2_\infty.}$$ Thus we obtain $$|||y||| \le 2C'|||F|||\le
4C'\|F\|_\infty  = 4C'\|F\|_{H^{\infty}_{\rm m}},\leqno(2.11)$$
similarly we have $({\bf E}\|y\|^2)^{1/2}\le C'\|F\|_2$ and by (2.9)
and (2.10) we conclude $$\|{\cal D}^{\eta}(F)\| \le
5C'K\|F\|_{H^{\infty}_{\rm m}}.\leqno(2.12)$$ This ends the proof of
Theorem 2.1.\qed Recapitulating, we have actually proved the following
statement.  \proclaim Theorem 2.2. Let $(\eta_n)$ be an adapted
sequence satisfying (2.2).  Let $$ {\cal D}^{\eta}:\ H^{\infty}_{\rm
m}\to B(H^2_{\rm m}(H),{H^2_{\rm m}(H)}^*)$$ be the mapping defined by
(2.3). Then ${\cal D}^{\eta}$ is a bounded linear operator satisfying
the identities (2.4) and (2.5). Moreover, for each $F$ in
$H^{\infty}_{\rm m}$, ${\cal D}^{\eta}(F)$  is ``Hankelian", meaning it
satisfies
 the following identity $$\forall G\in H^{\infty}_{\rm m}\quad {\cal
D}^{\eta}(F)M_{\rm m}(G)={}^tM_{\rm m}(G){\cal
D}^{\eta}(F).\leqno(2.13)$$

\pf Everything has already been explicitly proved, except (2.13) which
follows immediately from the definition (2.3) of ${\cal
D}^{\eta}(F)$.\qed \n {\bf Remarks.} The fact that the duality $H^1$,
BMO for martingales (see the classical reference [FS]) extends to the
Hilbert space valued case is a well known fact. The proof given in the
first pages of Garsia's book [Ga] extends almost verbatim. See also
 [Du, Pet].

Incidentally, the fact that the vectorial Nehari Theorem is valid in
this setting is also known (perhaps not so ``well''). All the
ingredients for this are available.  Indeed, it suffices to know that
any positive matrix valued function $W\colon \ \Omega \to M_n$ such
that $W\ge \vp I$ for some $\vp>0$ admits a factorization of the form
$W = F^*F$ with $F\colon \ \Omega\to M_n$ analytic and ``outer''. The
latter can be proved as in the case of the disc, see [HL].

\n In the disc case, the classical references for this are [Sa] and
[Pag].  Concerning vectorial Hankel operators, see  also [Pe3], and
consult [Tr] for more recent refinements. Finally, we also refer  to
[Pi2] for Banach space versions of Nehari's theorem.

\n {\bf Remark 2.3.} Although it is not needed to understand the
sequel, we should point out that  the arguments in the next section
are  based on properties of Brownian motion, following a well known
avenue in Probability Theory. Let $0\le r_0\le ...\le r_{n-1}<r_{n}<1$
be an increasing sequence of radii with $r_n\to 1$ when $n\to \infty$.

\n Let $(B_t)_{t>0}$ be the standard complex Brownian motion starting
from the origin and let $$T_n = \inf\left\{t>0\mid |B_t| = r_n\right\}
$$ be the usual stopping time.  We denote by $T=T_\infty$ the exit time
from the unit disc and by $({\cal A}_t)_{t>0}$ the Brownian
filtration.  Let ${\cal F}_n = {\cal A}_{T_n}$ for all $n\ge 0$.
Consider a function $F$ in $H^\infty$ and extend it analytically inside
$D$.  Then the martingale $(\psi_n)$ defined below is identical in
distribution with $B_{T_n}$ and $F(\psi_n)$ coincides in distribution
with $F(B_{T_n})$.
 Note that $B_{T_{n-1}}$ is uniformly distributed on $\{z\mid |z| =
r_{n-1}\}$.  Moreover for any fixed $z$ with $|z| = r_{n-1}$ the
distribution of $B_{T_n}$ conditional to $B_{T_{n-1}}=z$, is given by a
homothetic of the Poisson kernel. Thus, the formulae that we are using
below are entirely classical, they relate the Poisson kernels (and the
M\"obius transformations of the unit disc) with the distributions (or
certain conditional distributions) of the random variables $B_{T_n}$.
We refer the reader e.g. to [Du, Pet] or to the article [Bu] for more
information.

\vskip24pt \n {\bf \S 3. Proof of Theorem 1.1}

\def\m{{\rm m}} In this section, we deduce Theorem 1.1 from its
martingale counterpart Theorem 2.2.  the idea for this deduction is
quite simple: we introduce a specific element $\psi$ in $H^\infty_{\rm
m}$ with associated Hardy martingale $(\psi_n)_{n\ge 0}$ such that for
each $n\ge 1$ $$\forall Z\in \T\quad |\psi_n(Z)|=1-2^{-n}.$$ Then we
consider the unital homomorphism $$\sigma:\ A\to H^\infty_{\rm m}$$
defined, for $F$ in $A$, by the composition $$\sigma(F)=F(\psi).$$
Actually, we can also give a meaning to this definition for $F$ in
$H^{\infty}$: in that case we first extend $F$ inside the disc, then we
consider the uniformly bounded martingale $(F(\psi_n))$ and we define
$\sigma(F)$ as the weak-$*$ limit in $L^\infty(\Omega,P)$ of the
sequence $(F(\psi_n))$. When $F$ is in $A$, we recover $F(\psi)$ since
(by, say,  the martingale convergence theorem) $\psi_n$ tends to $\psi$
almost surely.

\n In addition, $\psi$ will be chosen uniformly distributed over the
circle, so that $\sigma$   extends by density to an isometric embedding
$$\sigma_2:\ H^2(H)\to H^2_\m(H).$$ Moreover, we will find an adapted
sequence $(\eta_n)_{n\ge 1}$ satisfying (2.2) and such that the
mappings $u$ and $u^\eta$, ${\cal D}^\eta$ defined respectively in
sections 1 and 2 satisfy for all polynomial $F$  in $A$ $$ {\cal
D}^\eta(\sigma(F))={\cal D}^\eta(\psi) M_\m ((\sigma(F')),\leqno(3.1)$$
$$  u(F)=  u^\eta(\sigma(F)).\leqno(3.2)$$ Then, by (2.5) and (2.12),
the mapping $\Gamma: H^2(H)\to H^2(H)^*$ defined by $$\Gamma  =
{}^t\sigma_2 {\cal D}^\eta(\psi) \sigma_2$$ satisfies all the
properties listed in Theorem 1.1.

 Consider first the disc at the origin. For any analytic function $F$
 in $A$ (considered as extended analytically inside $D$), we can write
$$F(0) = \int\limits_{\partial D} F(\xi) dm(\xi) \quad{\rm and}\quad
F'(0) = \int\limits_{\partial D} \bar\xi[ F(\xi)-F(0)] dm(\xi)$$ where
$dm(\xi)$ is normalized Lebesgue measure on $\partial D$.  The M\"obius
transformations  are usually denoted for any $z$ in $D$ by $\zeta\to
\ph_z(\zeta)$. For simplicity, we set
 $${\forall z\in D,\ \forall \zeta\in \T}\qquad
\Phi(z,\zeta)=\ph_z(\zeta)= {\zeta+z\over 1+\bar z
\zeta}.\qquad\qquad\qquad\qquad$$  Now consider $0\le r<s<1$ and $z$
with $|z|=r$. After a suitable change of variables, we can rewrite the
preceding formulae as follows:
  $$F(z)=\int F(s\Phi({z/s},\xi))dm(\xi) \leqno(3.3)$$ and
$$F'(z)s(1-{|z|^2\over {s^2}})=\int\bar\xi[F(s\Phi({z/s},\xi))-F(z)]
dm(\xi).  \leqno(3.4)$$

 \n For all $n\ge0$, let $$r_n=1-2^{-n},$$ so that $r_n<1$ and $r_n\to
1$.  We now define by induction a sequence of functions $(\psi_n)$ on
$\Omega=\T^I$.  We set first $\psi_0=0$ and, for all $n\ge 1$ and all
$Z=(Z_n)_{n\ge1}\in \Omega$
$$\psi_n(Z)=r_n\Phi(r_n^{-1}\psi_{n-1}(Z),Z_n).$$ Note that $\psi_n$
depends only on $Z_1,...,Z_n$, so that we can write (slightly
abusively) $\psi_n(Z)=\psi_n(Z_1,...,Z_n)$.  Note that for any $Z$ in
$\Omega$, the function $z\to \psi_n(Z_1,...,Z_{n-1},z)$ extends to an
analytic function in a neighborhood of $\bar D$ [namely $z\to
r_n\Phi(r_n^{-1}\psi_{n-1},z)$] such that
$\psi_n(Z_1,...,Z_{n-1},0)=\psi_{n-1}(Z_1,...,Z_{n-1})$.  Therefore
$(\psi_n)_{n\ge 0}$ is a Hardy martingale defining an element $\psi$ in
$H^\infty_{\rm m}$.  Moreover since $|\psi_n(Z)|\equiv r_n,$ we have
$\|\psi\|_{H^\infty_{\rm m}}=1$.  More generally, for any $F$ in $A$,
we have by (3.3) (=by Jensen's formula) $$\int
F(\psi_n(Z))dm(Z_n)=F(\psi_n(Z_1,...,Z_{n-1},0))=
F(\psi_{n-1}(Z)).\leqno(3.5)$$
Therefore, $(F\circ\psi_n)_{n\ge0}$ is a bounded Hardy martingale so
that $F\to \sigma(F)=F\circ \psi$ is a unital homomorphism from $A$ to
$H^\infty_{\rm m}$.  Consider $F$ as extended to an analytic function
on $D$. Integrating (3.5) we obtain $$\int F(\psi_n)dP =\int
F(\psi_{n-1})dP=\dots=F(0)$$ hence $$\int F(\psi_n)dP =\int
F(r_n\xi)dm(\xi).\leqno(3.6)$$ The latter identity ensures (since it
remains true for all harmonic functions in a neighborhood of $\bar D$)
that the random variable $\psi_n$ is uniformly distributed over the
circle $\{z\ | \ |z|=r_n\}$. Since $\psi_n$ tends a.s. to $\psi$,
$\psi$ is uniformly distributed over the unit circle.

\n Then (3.4) yields for all $n\ge 1$ $$F'(\psi_{n-1}(Z))
{{r_n^2-r^2_{n-1}}\over{r_n}}= \int
\bar{Z_n}[F(\psi_n(Z))-F(\psi_{n-1}(Z))] dm(Z_n).\leqno(3.7)$$ We have
clearly for all $r<1$ and for all $m\ge 0$ $$r^{m-1}m \widehat F(m)  =
\int e^{-i(m-1)t} F'(re^{it})dm(e^{it})$$ where $F = \sum\limits_{n\ge
0} \widehat F(n)z^n$ is an arbitrary function in $A$.  Taking $m=K_n$
and $r=r_{n-1}$ we obtain $$(r_{n-1})^{K_n-1} K_n \hat F(K_n) =\int
{\bar \xi}^{K_n-1}F'(r_{n-1}\xi)dm(\xi)$$ hence by (3.6) $$=\int
[r_{n-1}^{-1}\bar\psi_{n-1}]^{K_n-1} F'(\psi_{n-1}) dP.$$ Let
$\xi_{n-1}= [r_{n-1}^{-1} \psi_{n-1}]$. Note that $|\xi_{n-1}|=1$.
After multiplication of both sides of (3.7) by
$\bar{\xi}_{n-1}^{K_n-1}$ and integration with respect to $dm(Z_1)\dots
dm(Z_{n-1})$, we obtain for all $n>1$ $$\hat F(K_n)=\int\eta_{n-1}\bar{
Z}_n[F(\psi_n(Z))-F(\psi_{n-1}(Z))] dP(Z)\leqno(3.8)$$ with
$$\eta_{n-1}=\bar{\xi}_{n-1}^{K_n-1} . {{r_n}\over {r_n^2-r_{n-1}^2}}.
{{1}\over {K_n (r_{n-1})^{K_n-1}}}.$$ The first terms being irrelevant,
we may as well assume that $K_1=1$, so that taking $\eta_0=1$
identically, we still have (3.8) for $n=1$.  Since we assume
$2^{n-1}<K_n\le 2^n$ for $n>1$,  we clearly have (2.2) for some
numerical constant $C'$.

\n Now, since the distribution of $\psi_n$ is uniform over the circle
of radius $r_n$, we clearly  have, for all $F$ in $A$
$$\|F\|_{H^2}=\lim_{r\to1}\uparrow(\int |F(r\xi)|^2 dm(\xi))^{1/2}=
\sup_n\|F(\psi_n)\|_2=\|F\circ\psi\|_{H^2_{\rm m}},$$ so that
$\sigma_2$ as defined above is an isometry.  Moreover, (3.8)
immediately implies (3.2).

\n Finally, we verify (3.1). This is now easy. Indeed,  by (2.13) we
know that, for each $F$ in $A$, ${\cal D}(F)$ as defined above
satisfies (1.1) (\ie it is Hankelian), therefore by the remarks
preceding Theorem 1.1, we have necessarily (3.1) for any polynomial
$F$.
 This completes the proof of Theorem 1.1.\qed

\n {\bf Remark 3.1.} Let $F$  be an analytic function on $D$ with
boundary values
  in the classical space of BMO
 over the circle. Then it is well known that the martingale $F(\psi_n)$
is in the space $BMO_{\rm m}$  considered in the preceding section and
(with suitable choice of norms) we have
 $\|F(\psi)\|_{BMO_{\rm m}}\le \|F\|_{BMO}$. Therefore, we actually
have proved Theorem 1.1 with (0.5) replaced by the following stronger
inequality $$\exists~C~~~\forall~f\ {\rm polynomial}~~~~~~\N{\Gamma
M_{f'}}\leq C|||{f}|||\leqno(0.5)'$$ where (say) $$|||f||| =\sup_{z\in
D} (\int |f(\zeta)-f(z)|^2 P^z(d\zeta))^{1/2}.$$

\n {\bf Remark 3.2.} To prove Theorem 1.5, let us define $$\eta_{n-1,k}
=\bar{\xi}_{n-1}^{k} . {{r_n}\over {r_n^2-r_{n-1}^2}}.  {{1}\over {k
(r_{n-1})^{k-1}}}.$$ Note that there is a constant $C'$ such that
$\|\eta_{n-1,k}\|_{\infty}\le C'$ for all $n$ and all $k$ with
$2^{n-1}<k\le 2^n$. Moreover, we have for all $F$ in $A$ $$\hat
F(k)=\int\eta_{n-1,k}\bar{ Z}_n[F\psi_n(Z)-F\psi_{n-1}(Z)] dP(Z).$$ We
can now introduce the following modified version of ${\cal D}^{\eta}$:
for all
 $g,h$ in $H^2_{\rm m}(H)$ $${\cal D}(F)(g,h)  = {\bf E}
\left(\sum_{n\ge 1} \bar Z_n\, dF_n\, \sum_{2^{n-1}<k\le 2^n}m(k)
\eta_{n-1,k}\,  [C_kg,h] \right).$$ Then, the same argument as above
yields Theorem 1.5.

\n {\bf Remark 3.3.} There is a well known ``dictionary'' between the
classical theory of $H^1$ and BMO for the disc and the corresponding
theory for Brownian martingales. Although I did not see it, it seemed
very likely that there should be a way to prove that the operator
$\Gamma$ appearing in Theorem~1.1 satisfies (0.5) using only the
classical theory, and effectively this has recently been done by S.
Kisliakov (personal communication).

In the opposite direction, it is easy to modify the approach of \S 2 to
exhibit an ``Ito-Clifford'' integral (cf. [BSW])
 with the properties of Theorem~1.1.  \vskip24pt

\n {\bf \S 4. Operator space interpretations}

\def\m{{\rm m}} \def\I{{\cal I}} In this section, we give various
consequences of the previous results for operator spaces or (non
self-adjoint) operator algebras. In the theory of operator spaces, the
relevant morphisms are {\it completely bounded maps}   and the
corresponding isomorphisms are called  {\it complete isomorphisms}. See
e. g. [BP, ER,   Pa1 , Pi4].

 We start by an operator space theoretic reformulation of Theorem 1.1.
Fix a number $c>1$. Let $[A]_c$ be the  disc algebra equipped with the
operator space (actually operator algebra) structure induced by the
embedding $$A \subset \bigoplus_{\pi\in I(c)} B(H_\pi)$$ taking $a$ to
$\oplus_{\pi\in I(c)} \pi(a)$, where $I(c)$ is the class of all unital
homomorphisms $\pi\colon \ A\lgar B(H_\pi)$ with $\|\pi\|\le c$.  Note
that, as Banach algebras $[A]_c$ and $ A$ are the same with equivalent
norms, namely we have $\|a\|_A \le \|a\|_{[A]_c}\le c\|a\|_A$ for all
$a$ in $A$. However, as operator spaces,
 they are quite different. As we will see   this is a consequence of
Theorem 1.1 (and, by Theorem 0.1, it is  is equivalent to Corollary
1.2).

Indeed, consider the operator   space $\max(\ell_2)$ in the sense of
[BP].  The latter can be defined as follows.  Let $\I$ be the class of
all linear maps $v: \ell_2\to B(H_v)$ with $\|v\|\le 1$.  Let $J:
\ell_2\to \bigoplus_{v\in \I} B(H_v)$ be the isometric embedding
defined by $J(x)=\oplus_{v\in \I} v(x)$.  Then the operator space
$\max(\ell_2)$ can be defined as  the range of $J$.  Clearly, by
Theorem 1.1, the mapping $$P\colon \ [A]_c \lgar \max(\ell_2)$$ defined
by $PF = (\hat F(K_n))_{n\ge 1}$ is completely bounded (see the proof
of Theorem 4.1 below for details). Moreover, cf.\ e.g.\ [Pi1, p. 69]
$P$ is surjective.

Similarly, replacing $A$ by $A_\m$, we can define the operator algebra
$[A_m]_c$.
 Let us denote by $P_\m: A_\m \to \ell_2$ the   mapping defined by
$$\forall F\in A_\m \quad P_\m(F)= \sum_{n\ge 1} {\bf E}(dF_n\bar Z_n)
e_n.$$ In other words $P_\m$ is essentially the same as $u^\eta$ when
all $\eta_n$'s are identically equal to $1$. Then, it is easy to check
(by the same argument as above) that  $P_\m$ is a  bounded surjection
of $A_\m$ onto $\ell_2$. As observed in \S 2, it is not completely
bounded from $A_\m$ to $\max(\ell_2)$, but Theorem 2.2 implies that it
is completely bounded as a map from $[A_\m]_c$ to $\max(\ell_2)$.  More
generally, as an immediate consequence of Theorem 1.5, we have

\proclaim Theorem 4.1. Let $c>1$.  Let $(m(n))_{n\geq 0}$ be any scalar
sequence.  If  $(m(n))_{n\geq 0}$ satisfies (1.8), then the mapping
${\cal M}: [A]_c\to \max(\ell_2)$ defined (for $F\in A$) by ${\cal
M}(F)=(m(n)\hat F(n))$  is completely bounded.  However, if we view
${\cal M}$ as acting from $A$ to $\max(\ell_2)$, then it is completely
bounded iff $\sum|m(n)|^2<\infty$.  \def\I{\cal I}

\pf It suffices to show that for some constant $C$, we have for any
$v\in \I$, \break $\|v{\cal M}\|_{cb([A]_c,B(H_v))}\le C$.  Let
$C_n=v(e_n)$. Note that (1.5) holds, so that Theorem 1.1 and Remark 1.3
imply that there is a homomorphism $\pi\in I(c)$ and operators
$V_1,V_2$ so that  we can write,    for all $F$ in $A$ $$v{\cal
M}(F)=V_1\pi(F)V_2.$$ Moreover, we can achieve this with
$\|V_1\|\|V_2\|\le K_c$ for some constant $K_c$ depending only on $c$.
This implies $\|v{\cal M}\|_{cb([A]_c,B(H_v))}\le \|V_1\|\|V_2\|\le
K_c$, hence $\|{\cal M}\|_{cb([A]_c,\max(\ell_2))}\le K_c$.  This
proves the first part.
 The second one  follows from a joint  result of V. Paulsen and the
author (cf. [Pa3, Th. 4.1]).\qed

\proclaim Corollary 4.2. We have for any $c>1$  complete isomorphisms
$${[A]_c\over \hbox{Ker } P} \simeq \max(\ell_2) \quad {\rm and}\quad
 {[A_\m]_c\over \hbox{Ker } P_\m} \simeq \max(\ell_2).$$

\pf Indeed,  $[A]_c/\hbox{Ker }(P) \approx \ell_2$ as  Banach spaces,
so that the complete boundedness of the mapping $\max(\ell_2) \to
[A]_c/\hbox{Ker}(P)$ is obvious by maximality.  For the inverse
mapping,  the complete boundedness follows immediately from the
preceding statement. The case of $A_\m$ is analogous, we skip the
details.\qed

\n {\bf Remark.} The preceding two statements remain valid with
$H^\infty$ and $H^\infty_\m$ in the place of $A$ and $A_\m$.

\proclaim Corollary 4.3. When $c>1$, the operator
 spaces $[A]_c$ and $[A_\m]_c$ are not exact (in the sense of e.g.
[JP]).

\pf Indeed, by    [BP] we know that $\max(\ell_2)^*$ coincides with an
operator subspace of a commutative $C^*$-algebra (namely with
$\min(\ell_2)$). Hence, by [JP, Corollary 1.7], if $[A]_c$ was exact,
the mapping ${}^tP: \max(\ell_2)^* \to ([A]_c)^*$ would be $2$-summing,
which is absurd since it is an isomorphism on $\ell_2$. \qed

\n {\bf Remark.} Let $P^d_c$ be the subspace of $[A]_c$ spanned by the
polynomials of degree  at most $d$.  With the notation of [JP], the
preceding argument
 shows that for some absolute constant $\delta>0$ we have for all $c>1$
and all $d$ $$\delta {(c-1)}\sqrt{\ln(d)}\le d_{SK}(P^d_c).$$

In [Bo2], Bourgain proves an upper estimate for polynomially bounded
$n\times n$ matrices, which we will now be able to bound from below.
But first we consider another parameter (related to Theorem 0.1) which
we will
 estimate rather precisely. For any $n$, we will denote by $f(c,n)$ the
norm of the identity mapping from $M_n(A)$ to $M_n([A]_c)$.
Equivalently, we have $$f(c,n)=\sup \|(P_{ij}(T))\|_{M_n(B(H))}$$ where
the supremum runs over all polynomially bounded operators $T$ with
$\|u_T\|\le c$ and over all $n\times n$ matrices $(P_{ij})$ with
polynomial entries such that $\sup_{z\in\T}\|(P_{ij}(z))\|_{M_n}\le
1$.

\proclaim Theorem 4.4. There is an absolute constant $K>0$ such that
for all $c>1$ and all integers $n$ we have $$K^{-1} (c-1)\sqrt{n}\le
f(c,n)\le K c \sqrt{n}.\leqno(4.1)$$

\pf Let $(a_{ij})$ be an $n\times n$ matrix with entries in $B(H)$. It
is easy to check that $$\|{(a_{ij})}\|_{B(\ell_2^n(H))}\leq\sqrt
n\sup_i\|(\sum_j a_{ij} a_{ij}^*)^{1/2}\|_{B(H)}.\leqno(4.2)$$
   Let $(P_{ij})$ be an $n\times n$ matrix of polynomials as above such
   that $\sup_{z\in\T}\|(P_{ij}(z))\|_{M_n}\le 1$.  Then, for each $i$
and any $z$ in \T , we have $\sum_j|P_{ij}(z)|^2\leq 1$. By a result
due to Bourgain (see (21) in [Bo 2], this result uses [Bo1, Th. 2.2])
there is a numerical constant $K_1$ such that for each $i$ $$\|(\sum_j
P_{ij}(T)P_{ij}(T)^*)^{1/2}\|\leq K_1 c.$$ Hence by (4.2) we conclude
that $$f(c,n)\leq K_1 c\sqrt n.$$ To prove the converse direction in
(4.1), we will use the following known fact
  on random matrices (the idea to use this   in this context comes from
M. Junge's unpublished independent proof of the already mentioned
result from [Pa3]): there is a numerical constant $K_2$ such that for
each $n$ there is an $n$-tuple of unitary matrices $U_1,...,U_n$ in
$M_n$ such that $$\forall~(\alpha_i)\in\ell_2^n~~~~~~~\|{\sum\alpha_i
U_i}\|_{M_n}\leq K_2(\sum_1^n|\alpha_i|^2)^{1/2}.\leqno(4.3)$$ For a
proof see   e.g. [TJ, p. 323].  Let $C_i=(K_2)^{-1}U_i$ for
$i=1,2,...,n$ and $C_i=0$ for $i>n$. Then (1.5) is satisfied, but on
the other hand since the $U_i$'s are   (finite dimensional) unitaries,
we have $\|{\sum_1^n U_i\otimes\overline U_i}\|=n$ hence $\|{\sum_1^n
C_i\otimes \overline C_i}\|=(K_2)^{-2}n$.

\n Let $u: A\to M_n$ be the mapping associated to $(C_i)$ as in Theorem
1.1. Then by Theorem 1.1 and Remark 1.3,
  there is,  for some numerical constant $K_3$, a polynomially bounded
operator $T$ with $\|u_T\|\le c$ and operators $V_1,V_2$ with
$\N{V_1}\N{V_2}\leq K_3(c-1)^{-1}$ such that, for any polynomial $P$,
we have $$u(P)=V_1 P(T) V_2.$$ On the other hand, arguing as in the
proof of Corollary 1.2 we find that there is a matrix $(P_{ij})$ with
$\sup_{z\in\T}\|(P_{ij}(z))\|_{M_n}\le 1$, such that $$\N{(I\otimes
u)(P_{ij})}_{M_n(M_n)}\ge K_2^{-2}\sqrt {n}.$$ Hence this  implies
$$\N{(P_{ij}(T))}_{B(\ell_2^n(H))}\ge K_3^{-1}K_2^{-2}\sqrt n(c-1),$$
which yields the left side of (4.1) for a suitable constant $K$.\qed

In [B2], Bourgain proves that if $T\in M_N$ satisfies $\|u_T\|\le c$
there is an invertible $S\in M_N$ for which $\|S^{-1}TS\|\le 1$ and
satisfying $$\|S^{-1}\|\|S\| \le Kc^4 \ln(N+1)$$ for some absolute
constant $K$ (independent of $N$ or $c$).  By Theorem 0.1, this is
equivalent to $$\|u_T\|_{cb} \le Kc^4 \ln(N+1).$$ It is unclear how
sharp this estimate is asymptotically. However, as a simple direct
consequence of Theorem 1.1, we have

\proclaim Theorem 4.5. There is a constant $\delta>0$ with the
following property:  for any $N$ and $c>1$, there is $T\in M_N$
polynomially bounded with $\|u_T\|\le c$ such that any $S$ invertible
in $M_N$ with $\| S^{-1}TS\|\le 1$ must satisfy $$ \delta (c-1)
\sqrt{\ln(N+1)} \le \|S^{-1}\|\|S\|.$$

\pf Fix an integer $n$ and let $C_i$ be as in the preceding proof with
$\dim(H)=n$ and $C_i=0$ for all $i>n$. Let $\Gamma=\Gamma_\ph$ be as in
Theorem 1.1. For any polynomial $F$,   we again  let ${\cal
D}(F)=\Gamma_\ph M_{F'}$. Note that by Remark 1.6, we have
$\hat\ph(k)=0$ for all $k\le -2^n$.

\n Observe  that ${\cal D}(F)=\Gamma_\ph M_{F'}=\Gamma_{\ph{F'}}$, so
that ${\cal D}(F)=0$ if $F\in z^{2^n+1}A$.

\n Let us denote $K=H^2(H)/{z^{2^n+1}H^2(H)}$ (or equivalently
$K=H^2(H)\ominus {z^{2^n+1}H^2(H)}$) and for any $F$ in $A$ let $Q_F:
K\to K $ be the compression of $M_F$ to $K$.  The  observation
immediately preceding implies further that, for any polynomial $F$,
${\cal D}(F)$ vanishes
 on both $H^2(H)\times  {z^{2^n+1}H^2(H)}$ and $
{z^{2^n+1}H^2(H)}\times H^2(H)$. Therefore, ${\cal D}(F)$ defines
unambiguously a linear map $\Delta(F): K\to K^*$ satisfying
$\|\Delta(F)\| \le \|{\cal D}(F)\|$ hence (by Theorem 1.1) we have for
any polynomial $F$ $$\|\Delta(F)\| \le C \|F\|_\infty.$$ Finally, let
${\cal H}=K^*\oplus K$, let $\eps=C^{-1}(c-1)$ and let $r(F): {\cal
H}\to {\cal H}$ be defined by
$$r(F)=\left(\matrix{{^tQ_F}&{\eps\Delta}(F)\cr 0& Q_F\cr}\right)$$
Then, it is easy to check that $F\to r(F)$ is a bounded homomorphism
with $\|r\|\le c$ for which there are contractions $V_1$ and $V_2$ such
that $V_1 r(z^{2^k}) V_2=\eps C_k$ forall $k=1,2,...,n$. By the same
argument as in the preceding proof, we know that this implies
$$\|r\|_{cb} \ge \eps (K_2)^{-2} \sqrt{n}.\leqno(4.4)$$ But on the
other hand we have $\dim({\cal
H})=2\dim({K})=2(2^n+1)\dim({H})=2(2^n+1)n$ so that $n\approx
\ln(\dim({\cal H})$, which shows that (4.4) yields the announced
result, modulo Paulsen's criterion (cf. Theorem 0.1). \qed

\n {\bf Acknowledgement.} I am very grateful to Quanhua Xu who,
 in answer to a question of mine, showed me the above proof of (2.11),
at a time when, perhaps, I was being
 paralyzed by pessimism.  At that stage, knowing (2.11), I could almost
immediately complete
 the original version of the present investigation, which, by the way,
did not directly involve Hankel operators.  (Shortly after, I realized
that the example could be described as a vectorial Hankel operator, as
above.) I am also grateful to Xu for his help in checking the first
drafts. I also thank S. Treil and B. Maurey for  pertinent remarks on
intermediate versions.

\bigskip\bigskip\bigskip \centerline {\bf REFERENCES} \bigskip\bigskip

\item{AP} A. B. Aleksandrov and V. V. Peller. Hankel operators and
similarity to a contraction. Preprint, 1995.

 \item{BR} O. Bratelli and D. Robinson. Operator algebras and quantum
statistical mechanics II. Springer Verlag, New-York, 1981.

\item{BSW} C.\ Barnett, R.F.\ Streater, I.F.\ Wilde, The It\^o-Clifford
Integral, J.\ Funct.\ Analysis 48 (1982), 172--212; \item{} The
It\^o-Clifford Integral II - Stochastic differential equation,
J.\ London Math.\ Soc. 27 (1983), 373--384; \item{} The It\^o-Clifford
Integral III -- Markov property of solutions to Stochastic differential
equation, Commun.\ Math.\ Phys. 89 (1983), 13--17; \item{} The
It\^o-Clifford Integral IV -- A Radon-Nikodym theorem and bracket
processes, J.\ Operator Theory 11 (1984), 255--211.

\item{Bo1} J. Bourgain. New Banach space properties of the disc algebra
and $H^\infty$. Acta Math. 152 (1984) 1-48.

\item{Bo2}  J.  Bourgain.  On the similarity problem for polynomially
bounded operators on Hilbert space.  Israel J.  Math.  54 (1986)
227-241.

 \item{BP} D. Blecher and V. Paulsen. Tensor products of operator
 spaces.  J. Funct. Anal. 99 (1991) 262-292.

\item{Bu} D. Burkholder. Distribution function inequalities for
martingales.  Ann. Probab. 1 (1973) 19-42.

\item {D}   P. Duren. The theory of $H^p$ spaces. Academic Press,
New-York, 1970.

 \item{Du} R. Durrett.  Brownian motion and martingales in analysis.
Wadsworth Math. Series, Belmont (California) 1984.

\item{ER} E. Effros and Z.J. Ruan.  A new approach to operator spaces.
 Canadian Math. Bull.  34 (1991) 329-337.

\item{Fo}  S. Foguel.   A counterexample to a problem of Sz. Nagy.
Proc. Amer. Math. Soc.  15 (1964) 788-790.

\item{FS} C. Fefferman and E. Stein. $H^p$-spaces in several variables.
Acta Math. 129 (1972) 137-193.

\item{FW}  C. Foias  and J. P.   Williams.    On a class of
polynomially bounded operators.   Preprint  (unpublished, 1979 or 1980
?).

\item{Ga} A. Garsia.  Martingale inequalities: seminar notes on recent
progress.  Benjamin, 1973.

\item{H} U. Haagerup. Injectivity and decomposition of completely
 bounded maps in ``Operator algebras and their connection with Topology
 and Ergodic Theory''. Springer Lecture Notes in Math. 1132 (1985)
170-222.

\item{Ha1}    P. Halmos.  Ten problems in Hilbert space.   Bull.  Amer.
Math. Soc.        76 (1970) 887-933.

\item{Ha2}    P. Halmos.  On Foguel's answer to Nagy's question. Proc.
Amer. Math. Soc.  15 (1964) 791-793.

\item{HL} H. Helson  and H. Lowdenslager.  Prediction theory and
Fourier series in several variables.
 Acta Math  99 (1958) 165-202.

\item{JP} M. Junge and G. Pisier. Bilinear
 forms on exact operator spaces and $B(H)\otimes B(H)$. Geometric
and Functional Analysis (GAFA Journal) 5 (1995) 329-363.

\item{Le}   A. Lebow.   A power bounded operator which is not
polynomially bounded.   Mich. Math. J.  15 (1968) 397-399.

\item{LPP}  F. Lust-Piquard  and G.  Pisier.   Non commutative
Khintchine and Paley inequalities.   Arkiv f=9Ar Mat.   29 (1991)
241-260.

\item{vN}  J.  von Neumann.   Eine spektraltheorie f\"ur allgemeine
operatoren eines unit\"aren \break raumes. Math. Nachr.  4 (1951)
49-131.

\item{Ni}  N.  Nikolskii.    Treatise on the shift operator.  Springer
Verlag, Berlin 1986.

\item{Pa1}  V.  Paulsen.    Completely bounded maps and dilations.
Pitman Research Notes in Math. 146, Longman, Wiley, New York, 1986.

\item{Pa2}  V. Paulsen.    Every completely polynomially bounded
operator is similar to a contraction.   J. Funct. Anal.  55 (1984)
1-17.

\item{Pa3}  V. Paulsen.   The maximal operator space of a normed
space.
 Proc. Edinburgh Math. Soc. To appear.

 \item{Pag}   L. Page.   Bounded and compact vectorial Hankel
operators.   Trans. Amer. Math. Soc.   150 (1970) 529-540.

\item{Pe1}  V.  Peller.   Estimates of functions of power bounded
operators on Hilbert space J. Oper. Theory   7 (1982) 341-372.

\item{Pe2}   V. Peller.  Estimates of functions of Hilbert space
operators, similarity to a contraction and related function algebras.
in Linear and complex analysis problem book (edited by Havin and
Nikolskii) Springer Lecture notes 1573 (1994) 298-302.

\item{Pe3}   V. Peller. Vectorial Hankel operators, commutators and
related operators of the Schatten von Neumann classes. Integral
Equations and Operator Theory 5 (1982) 244-272.

\item{Pet} K. Petersen. Brownian motion, Hardy spaces and bounded mean
oscillation.  LMS Lecture notes series 28, Cambridge Univ. Press,
Cambridge, 1979.

\item{Pi1}   G.   Pisier. Factorization of linear operators and the
Geometry of Banach spaces.  CBMS (Regional conferences of the A.M.S.)
no 60, (1986) Reprinted with corrections 1987.

\item{Pi2}  G.  Pisier. Factorization of operator valued analytic
functions.  Advances in Math.   93 (1992) 61-125.

\item{Pi3}    G. Pisier.        Multipliers and lacunary sets in
non-amenable groups.   Amer. J. Math.   117 (1995) 337-376.

\item{Pi4}   G. Pisier. Similarity problems and completely bounded
maps.  Springer Lecture notes 1618 (1995).

\item{Ro} R. Rochberg.   A Hankel type operator arising in deformation
theory.  Proc. Symp. Pure Math.  35 (1979)  457-458.

\item{Ru}   W.   Rudin. Fourier analysis on groups.  Interscience. New
York, 1962.

\item{Sa}  D.    Sarason. Generalized interpolation in $H^{\infty}$.
Trans. Amer. Math. Soc.  127 (1967) 179-203.

\item{SN}   B. Sz.-Nagy. Completely continuous operators with uniformly
bounded iterates.  Publ. Math. Inst.  Hungarian Acad. Sci.  4 (1959)
89-92.

\item{SNF}  B. Sz.-Nagy  and  C.        Foias. Harmonic analysis of
operators on Hilbert space   Akademiai Kiad\'o, Budapest   1970.

\item{St}  J. Stafney. A class of operators and similarity to
contractions.  Michigan Math. J. 41 (1994) 509-521.

\item{TJ} N. Tomczak-Jaegermann.  Banach-Mazur distances  and finite
dimensional operator ideals. Longman, Pitman monographs and surveys in
pure and applied math. 38 (Wiley) 1989.

\item{Tr} S. Treil. Geometric methods in spectral theory of vector
valued functions: some recent results. in: Operator Theory: Adv. Appl.
42 (1989) 209-280. (Birkhauser)

\vskip12pt

Texas A\&M University, College Station, TX 77843, U. S. A.

and

Universit\'e Paris VI, Equipe d'Analyse, Case 186, 75252 Paris Cedex
05, France

  {\bf GIP@CCR.JUSSIEU.FR}

\end \bye